\documentclass[a4paper,11pt]{amsart}
\usepackage{amsfonts,amssymb,amsthm,cite,amsmath,amstext}
\usepackage{color}
\usepackage[colorlinks,linkcolor=black,citecolor=black]{hyperref}
\usepackage{comment}
\usepackage{framed}
\usepackage{calligra} 
\usepackage{tikz-cd}

\definecolor{shadecolor}{gray}{0.875}

\newtheorem{thrm}{Theorem}[section]
\newtheorem{thrmx}{Theorem}

\newtheorem{lem}[thrm]{Lemma}
\newtheorem{cor}[thrm]{Corollary}

\newtheorem{conj}[thrm]{Conjecture}

\theoremstyle{definition}

\newtheorem{exmple}[thrm]{Example}
\newtheorem{rmk}[thrm]{Remark}

\newenvironment{claim}
            {\par \bigskip \noindent \textbf{Claim:}}
            {$\Box$ \par \noindent}

\title{The inequalities of Chern classes and Riemann-Roch type inequalities}
\author{Xing Lu and Jian Xiao}
\date{}

\begin{document}
\maketitle


\begin{abstract}
Motivated by Koll\'{a}r-Matsusaka's Riemann-Roch type inequalities, applying effective very ampleness of adjoint bundles on Fujita conjecture and log-concavity given by Khovanskii-Teissier inequalities, we show that for any partition $\lambda$ of the positive integer $d$ there exists a universal bivariate polynomial $Q_\lambda(x, y)$ which has $\deg Q_\lambda \leq d$ and whose coefficients depend only on $n$ and $\lambda$, such that for any projective manifold $X$ of dimension $n$ and any ample line bundle $L$ on $X$,
\begin{equation*}
  \left|c_\lambda(X)\cdot L^{n -d}\right|\leq
  \frac{Q_{\lambda}(L^{n}, K_X \cdot L^{n -1} )}{(L^{n})^{d-1}},
\end{equation*}
where $K_X$ is the canonical bundle of $X$ and $c_\lambda(X)$ is the monomial Chern class given by the partition $\lambda$. As a special case, when $K_X$ or $-K_X$ is ample, this implies that there exists a constant $c_n$ depending only on $n$ such that for any monomial Chern classes of top degree, the Chern number ratios satisfy the following inequality
\begin{equation*}
\left|\frac{c_\lambda(X)}{c_1 (X) ^{n}}\right|\leq c_n,
\end{equation*}
which recovers a recent result of Du-Sun. The main result also yields an asymptotic version of the sharper Riemann-Roch type inequality. Furthermore, using similar method we also obtain inequalities for Chern classes of the logarithmic tangent bundle.
\end{abstract}

\tableofcontents

\section{Introduction}
Throughout this paper, we work on the field of complex numbers $\mathbb{C}$.

We always assume that $X$ is a projective manifold of dimension $n$, and that $L$ is an ample line bundle on $X$. The canonical bundle of $X$ is denoted by $K_X$.

We denote by $P(d,n)$ the set of all partitions of the positive integer $d$ by nonnegative integers no more than $n$:
\begin{equation*}
  P(d,n)=\{\lambda = (\lambda_{1},...,\lambda_{r}) \in \mathbb{N} ^r: |\lambda|=\sum_{i=1}^r \lambda_i = d, \lambda_i \leq n, \ \forall 1\leq i \leq r\}.
\end{equation*}
In the sequel, we always assume $d\leq n$. Given a partition $\lambda \in P(d,n)$, we denote
$$c_{\lambda}(X)  = \prod_{i=1}^{r} c_{\lambda_{i}}(X) \in H^{2d}(X,\mathbb{Z}),$$
where $c_{\lambda_{i}}(X)$ is the $\lambda_{i}$-th Chern class of $X$. We call $c_{\lambda}(X)$ the monomial Chern class corresponding to $\lambda$. In the intersection numbers involving $c_1 (M)$ of a line bundle $M$, we simply denote $c_1 (M)$ by $M$.

\subsection{Motivation}
A fundamental problem (Riemann-Roch problem) in algebraic geometry is to find the exact value of $h^0 (X, kL)$ of a line bundle $L$ on a projective manifold $X$ of dimension $n$. When $L$ is ample and $k$ is large enough, by Kodaira vanishing theorem and Hirzebruch-Riemann-Roch theorem, the value is given by the Hilbert polynomial
\begin{align*}
  h^0 (X, kL)&= \chi(X,kL) = \int_X ch(kL) \cdot td(X)\\
  &=\frac{L^n}{n!} k^n -\frac{K_X \cdot L^{n-1}}{2(n-1)!} k^{n-1}+...,
\end{align*}
whose coefficients are determined by the intersection numbers of the monomial Chern classes and $L$:
$$ c_{\lambda}(X) \cdot L^{n-d} : \lambda \in P(d,n), 1 \le d \le n .$$
When the value of $k$ is arbitrary, a fundamental result by Koll\'ar and Matsusaka \cite{kollarMatRiemRoch} implies that there exists a universal polynomial $Q(z)$ with $\deg_{z}Q \le n-1$, such that for all $k>0$,
$$\left| h^{0}(X,kL)-\frac{L^{n}}{n!}k^{n} \right | \le Q(k),$$
whose coefficients are determined by $L^{n},K_{X} \cdot L^{n-1}, n$.

Koll\'ar-Matsusaka's theorem had been generalized by Luo \cite{luoRiemRoch} to the case when $X$ is a smooth projective variety and $L$ is big and nef.
Indeed, Koll\'ar-Matsusaka's result also holds whenever $X$ is a normal projective variety and $L$ is big and semiample. It was refined by Nielsen \cite{nielRRineq} when $X$ is a normal projective variety and $L$ is ample, by proving that there exists a universal polynomial $Q(z)$ with $\deg Q \le n-2$, such that for all $k>0$,
$$\left| h^{0}(X,kL)-\frac{L^{n}}{n!}k^{n}+ \frac{K_X \cdot L^{n-1}}{2(n-1)!} k^{n-1}\right | \le Q(k).$$
The coefficients of $Q(z)$ are determined by $L^{n},K_{X} \cdot L^{n-1}, n$.

In the ample case, combining the Riemann-Roch inequality and the Hirzebruch-Riemann-Roch theorem, it is natural to ask:
\begin{quote}
\emph{Given $\lambda \in P(d,n)$, is the Chern number $c_{\lambda}(X) \cdot L^{n-d}$ effectively and explicitly controlled by the values of $L^{n}$ and $K_{X} \cdot L^{n-1}$?}
\end{quote}

This can be seen as a somehow stronger version of the asymptotic Riemann-Roch inequality.
The aim of this note is to give a confirmative answer to this question.

\subsection{Main results}
We prove the following result.

\begin{thrmx}\label{main thrm1}
For any partition $\lambda \in P(d, n)$, there exists a universal bivariate polynomial $Q_\lambda(x, y)$ satisfying that
\begin{itemize}
  \item the degree of $Q_\lambda$, $\deg Q_\lambda \leq d$,
  \item the coefficients of $Q_\lambda$ depend only on $n$ and $\lambda$,
\end{itemize}
such that for any projective manifold $X$ of dimension $n$ and any ample line bundle $L$ on $X$,
\begin{equation*}
  \left|c_\lambda(X)\cdot L^{n -d}\right|\leq
  \frac{Q_{\lambda}(L^{n}, K_X \cdot L^{n -1} )}{(L^{n})^{d-1}}.
\end{equation*}
In particular, $c_\lambda(X)\cdot L^{n -d}$ is polynomially controlled by $L^{n}, K_X \cdot L^{n -1}$ in the sense that
$$\left|c_\lambda(X)\cdot L^{n -d}\right|\leq
 Q_{\lambda}(L^{n}, K_X \cdot L^{n -1} ).$$
\end{thrmx}

From the proof, one can indeed obtain an expression of the polynomial $Q_\lambda (x,y)$. The bound given by $Q_\lambda (x,y)$ can be explicit, but it may be rough in general.
The dependence on $L^{n}$ and $K_X \cdot L^{n -1}$ is similar to the effective Matsusaka's big theorem \cite{matsBigThrm, siuEffMatsAIF, demaiEffectiveVeryAmple}.

We establish Theorem \ref{main thrm1} by the following more finer estimates.

\begin{thrmx}\label{main thrm2}
Fix a positive integer $n$.
\begin{enumerate}
  \item For any partition $\lambda \in P(d, n)$, there exist linear polynomials $P_{\lambda}^{\pm}(x_{0},...,x_{d})$ whose coefficients depend only on $n$ and $\lambda$, such that for any projective manifold $X$ of dimension $n$ and any ample line bundle $L$ on $X$,
      \begin{equation*}
        P_{\lambda}^{-}(K_{X}^{i} \cdot L^{n-i} : 0 \le i \le d) \le c_{\lambda}(X) \cdot L^{n-d} \le P_{\lambda}^{+}(K_{X}^{i} \cdot L^{n-i} : 0 \le i \le d).
      \end{equation*}

  \item For any $2 \le i \le n$, there exist polynomials $R_{i}^{\pm}(x,y)$ whose degrees $\deg R_{i}^{\pm} \le i$ and whose coefficients depend only on $n$ and $i$, such that for any projective manifold $X$ of dimension $n$ and any ample line bundle $L$ on $X$,
      \begin{equation*}
        R_{i}^{-}(L^{n},K_{X} \cdot L^{n-1}) \le (K_{X}^{i} \cdot L^{n-i}) (L^{n})^{i-1} \le R_{i}^{+}(L^{n},K_{X} \cdot L^{n-1}).
      \end{equation*}
\end{enumerate}

\end{thrmx}

\subsection{Applications and extensions}

\subsubsection{Boundedness for Chern numbers}
In the recent work by Du-Sun \cite{durongChernIneq}, by using the method of pulling back Schubert classes in the Chow group of a Grassmannian under the Gauss map, it was proved that: for a projective manifold $X$ over an algebraically closed field $k$, if $K_X$ or $-K_X$ is ample and $k$ is of characteristic 0, or if $K_X$ or $-K_X$ is ample and globally generated and $k$ is of positive characteristic, then there exists a constant $c_n$ depending only on $n$ such that for any monomial Chern classes of top degree, the Chern number ratios satisfy the following inequality 
\begin{equation*}
\left|\frac{c_\lambda(X)}{c_1 (X) ^{n}}\right|\leq c_n.
\end{equation*}

The result of Du-Sun gave an affirmative answer to a generalized open question, that whether the region described by the Chern ratios is bounded, posted by Hunt \cite{huntchern}. Theorem \ref{main thrm1} generalizes their result over an algebraically closed field of characteristic 0, in the sense that one can drop the assumption of being Fano or of general type and it holds on any projective manifolds by introducing the ``twist'' ample line bundle $L$. In particular, letting $L=K_X$ or $L=-K_X$ immediately recovers \cite{durongChernIneq}. This provides an alternative approach to Du-Sun's boundedness theorem for Chern ratios.

\begin{cor}\label{chern ratio}
Let $X$ be a projective manifold of dimension $n$, with $K_{X}$ or $-K_{X}$ ample, then the Chern number ratios
$$\left [ \frac{c_{\lambda}(X)}{c_{1}(X)^{n}} \right ]_{\lambda \in P(n,n)} \in \mathbb{R}^{p(n)}$$
is contained in a bounded set in $\mathbb{R}^{p(n)}$ independent of $X$, where $p(n)$ is the partition number.
\end{cor}

\begin{rmk}
In the Fano case, as pointed out to us by Haidong Liu, the above boundedness can be derived from the boundedness of Fano manifolds \cite{KMM92}. By \cite{birkarBAB}, it also extends to certain Chern classes on singular Fano varieties. In particular, by \cite{liuiwai2023miyaoka} for any terminal weak Fano variety of dimension $n$, there is a constant $b_n$ depending only on $n$ such that
$$c_1 (X)^n \leq b_n c_2 (X) \cdot c_1 (X)^{n-2}. $$
For more refined estimates regarding this aspect, we refer the reader to the recent papers \cite{liu2023kawamatamiyaoka}, \cite{liuiwai2023miyaoka}, \cite{liujieC2Fano} and the references therein.
\end{rmk}

Regarding the boundedness, we have:
\begin{cor}
There is some uniform constant $c(n, v, w) >0$ such that for any projective manifold $X$ of dimension $n$ and any ample line bundle $L$ on $X$ with
\begin{itemize}
  \item $L^n \leq v$,
  \item $K_X \cdot L^{n-1} \leq w$,
\end{itemize}
we have
\begin{equation*}
  \left|c_\lambda(X)\cdot L^{n -d}\right|\leq c(n, v, w).
\end{equation*}

\end{cor}

This is compatible with the boundedness result by \cite[Theorem 1.4]{kollarEffBasepiont} (see also \cite[Theorem 3]{kollarMatRiemRoch}). The point from our approach is that the bound $c(n, v, d)$ can be theoretically explicit, given by the polynomial in Theorem \ref{main thrm1}.

\subsubsection{Asymptotic Riemann-Roch type inequalities}
Another consequence is an asymptotic version of the sharper Riemann-Roch type inequality for ample line bundles.

\begin{cor}\label{sharper rr}
There exists a polynomial $Q(z)$ with $\deg Q \le n-m-1$, such that for any projective manifold $X$ of dimension $n$ and any ample line bundle $L$ on $X$, for $k$ large enough,
$$\left| h^{0}(X,kL)-\sum_{i=0}^m a_i k^{n-i}\right | \le Q(k),$$
where $\sum_{i=0}^m a_i k^{n-i}$ is the truncation of the Hilbert polynomial
$$ \chi(X,kL) = \int_X ch(kL) \cdot td(X)=\sum_{i=0}^n a_i k^{n-i}.$$
The coefficients of $Q(z)$ are determined by $L^{n},K_{X} \cdot L^{n-1}, n$ and the partitions $\lambda$ with $|\lambda|\geq m+1$.
\end{cor}

It improves \cite[Lemma 5.2]{kollarMatRiemRoch}, where Koll\'ar-Matsusaka proved that every $|a_i|$ is bounded above by a polynomial of $L^{n},K_{X} \cdot L^{n-1}$ of degree $n$. When $m=1$, this gives the asymptotic version of the sharper Riemann-Roch type inequality obtained by Nielsen \cite{nielRRineq}.

\subsubsection{Logarithmic tangent bundles}
While not a straightforward application,
using similar method, we extend Theorem \ref{main thrm1} to the logarithmic case (with less precise estimates).
We prove the following result:

\begin{thrm}\label{main thrm log}
Let $X$ be a projective manifold of dimension $n$, and let $D=\sum_{i=1} ^N D_i$ be a reduced SNC divisor on $X$. Denote by $\Omega_X ^1 (\log D)$ the logarithmic cotangent bundle of the pair $(X, D)$. Let $L$ be an ample line bundle on $X$. Given $J=(j_1,...,j_N)\in \mathbb{N}^N$, denote $D^J = D_1 ^{j_1}\cdot...\cdot  D_N ^{j_N}$ and $|J|= j_1 +...+j_N$. Then for any partition $\lambda \in P(d, n)$, the Chern number $c_\lambda (\Omega_X ^1 (\log D))\cdot L^{n-d}$ is polynomially controlled by the intersection numbers $$x_J=D^J\cdot L^{n-|J|},\  y_J=K_X \cdot D^J\cdot L^{n-|J|-1},$$
where $J$ varies among the index set $\{J= (j_1,...,j_N): 0\leq j_k \leq 2\}$. More precisely, there is a positive integer $k(\lambda, n)$ and
a family of polynomials $\{Q_{\lambda} ^l\}_{l=1} ^{k(\lambda, n)}$ such that
\begin{align*}
  \min_{1\leq l\leq k(\lambda, n)} \{Q_{\lambda} ^l(x_J, y_J)\}
  \leq c_\lambda (\Omega_X ^1 (\log D))\cdot L^{n-d}
  \leq \max_{1\leq l\leq k(\lambda, n)} \{Q_{\lambda} ^l (x_J, y_J)\}.
\end{align*}

\end{thrm}

In particular, using the fact that $c_n (\Omega_X ^1 (\log D)) = (-1)^n \chi (X\setminus D)$, where $\chi (X\setminus D)$ is the Euler characteristic of $X\setminus D$, Theorem \ref{main thrm log} gives an intersection theoretic bound for the Euler characteristic of a smooth quasi-projective variety by the boundary and the canonical divisor of its SNC compactification.

This paper is organized as follows. In Section \ref{prel} we present the key ingredients that shall be applied in the proof. Section \ref{sec main} is devoted to the proof of the main theorems and its applications.

\subsection*{Acknowledgements}
This work is supported by the National Key Research and Development Program of China (No. 2021YFA1002300) and National Natural Science Foundation of China (No.
11901336). We would like to thank Haidong Liu for helpful comments. We also thank the referee for the careful reading and helpful comments.

\section{Preliminaries}\label{prel}

For the basics on the positivity in algebraic geometry, we refer the reader to \cite{lazarsfeldPosI, lazarsfeldPosII}.

In our approach to the inequalities of Chern classes, we apply the geometry of adjoint bundles in an essential manner.
This is closely related to the tantalizing Fujita conjecture:
\begin{conj}
Let $X$ be a projective manifold of dimension $n$ and $L$ an ample line bundle on $X$, then
\begin{itemize}
  \item Freeness conjecture: for $m\geq n+1$, $K_X + (n+1)L$ is free.
  \item Very ampleness conjecture: for $m\geq n+2$, $K_X + (n+2)L$ is very ample.
\end{itemize}
\end{conj}

There are a lot of important works towards this conjecture, the works that we mention below is far from complete.
When $n=2$, the conjecture was completely settled by Reider \cite{reidFujita}. For general dimension, the first breakthrough was made by Demailly in the seminal work \cite{Dem93}, who proved that $2K_X + mL$ is very ample whenever $m\geq 12n^n$ (see also \cite{einlazarZeroEstimates} for the algebro-geometric analogues of some of the
analytic parts of Demailly's argument).
This result has also been reproved and somehow improved in \cite{siuVeryAmple, demaiEffectiveVeryAmple}, in particular, $2K_X + mL$ is very ample whenever $m\geq 2 + \binom{3n+1}{n}$. For the freeness part, when $n\leq 5$ it was solved in a series of works \cite{einlazFree3, kawamaFree34, yezhuFree5}; for general dimension, the major achievement was made by Angehrn-Siu \cite{siuEffFree} who proved that $K_X + m L$ is free whenever $m\geq \frac{n^2 +n +2}{2}$. On the other hand, by Mori's theory and the base point theorem \cite{kollarmoriBOOK}, $K_X + (n+1)L$ is always semiample. In the case when $L$ has stronger positivity, as examples we refer the reader to the classical works \cite{karenfujitaLocal, karenfujita2} that work in any characteristic.

We summarize the following fact that is enough in our setting.

\begin{lem}\label{veryampleness}
Let $X$ be a projective manifold of dimension $n$ and $L$ an ample line bundle on $X$, then \begin{itemize}
  \item The adjoint bundle $K_X + (n+1)L$ is nef. In particular, $K_X + (n+2)L$ is ample.
  \item There is some explicit $C_n$ such that $2K_X + mL$ is very ample whenever $m\geq C_n$. Explicitly, one can take $C_n =2 + \binom{3n+1}{n}$.
\end{itemize}

\end{lem}

In our setting, we shall apply the positivity of twisted tangent bundle. The following simple lemma is also quite useful.

\begin{lem}[see Lemma 12.1 of \cite{Dem93}] \label{tangent twist}
Let $X$ be a projective manifold of dimension $n$ and $A$ a very ample line bundle on $X$. Then the vector bundle $T_{X} \otimes \mathcal{O}_{X}(K_{X} + nA)$ is globally generated, and hence nef.
\end{lem}

\begin{exmple}
The above result is optimal in the following sense.
Let $X = \mathbb{P}^{n}$ with $n\geq 2$ and $A=\mathcal{O}(1)$ the hyperplane line bundle on $\mathbb{P}^{n}$, then for any $\epsilon>0$,
$$E=T_{\mathbb{P}^{n}} \otimes (K_{\mathbb{P}^{n}} + (n-\epsilon)\mathcal{O}(1))$$
is not nef. To check this, one just need to restrict $E$ to some line in $\mathbb{P}^{n}$ and then apply Grothendieck's splitting theorem.
\end{exmple}

Regarding the positivity of Chern classes of nef vector bundles, we have:

\begin{lem}[see \cite{lazarsfeldPosII}, or Corollary 2.6 of \cite{DPSneftangent}] \label{nef chernclass}
Let $X$ be a projective manifold of dimension $n$ and $E$ a nef vector bundle on $X$. Then for any monomial Chern class $c_\lambda (E)$ of degree $2k$ and ample classes $L_1,...,L_{n-k}$,
$$c_{1}(E)^{k} \cdot L_1\cdot...\cdot L_{n-k} \ge c_{\lambda}(E)\cdot L_1\cdot...\cdot L_{n-k} \ge 0.$$
\end{lem}

\begin{rmk}
In the K\"ahler setting, by \cite{DPSneftangent} Lemma \ref{nef chernclass} also holds when $L_1 =...=L_{n-k} =\omega $ is a K\"ahler class.
\end{rmk}

Another result that we will apply is the Khovanskii-Teissier inequality (see e.g. \cite[Section 1.6]{lazarsfeldPosI}):

\begin{lem} \label{ktineq}
Let $X$ be a projective manifold of dimension $n$ and let $A, B, C_1,...,C_{n-m}$ be nef divisor classes on $X$, then the sequence $$\{A^k \cdot B^{m-k}\cdot C_1 \cdot...\cdot C_{n-m}\}_{k=0} ^m$$ is log-concave.
\end{lem}


\section{Proof of the main results}\label{sec main}

In this section, we give the proof of the main results.

First, we prove:

\begin{thrm}\label{thrm linear}
Fix a positive integer $n$.
\begin{enumerate}
  \item For any partition $\lambda \in P(d, n)$, there exist linear polynomials $P_{\lambda}^{\pm}(x_{0},...,x_{d})$ whose coefficients depend only on $n$ and $\lambda$, such that for any projective manifold $X$ of dimension $n$ and any ample line bundle $L$ on $X$,
      \begin{equation*}
        P_{\lambda}^{-}(K_{X}^{i} \cdot L^{n-i} : 0 \le i \le d) \le c_{\lambda}(X) \cdot L^{n-d} \le P_{\lambda}^{+}(K_{X}^{i} \cdot L^{n-i} : 0 \le i \le d).
      \end{equation*}

  \item For any $2 \le i \le n$, there exist polynomials $R_{i}^{\pm}(x,y)$ whose degrees $\deg R_{i}^{\pm} \le i$ and whose coefficients depend only on $n$ and $i$, such that for any projective manifold $X$ of dimension $n$ and any ample line bundle $L$ on $X$,
      \begin{equation*}
        R_{i}^{-}(L^{n},K_{X} \cdot L^{n-1}) \le (K_{X}^{i} \cdot L^{n-i}) (L^{n})^{i-1} \le R_{i}^{+}(L^{n},K_{X} \cdot L^{n-1}).
      \end{equation*}
\end{enumerate}

\end{thrm}

By Lemma \ref{veryampleness} and Lemma \ref{nef chernclass}, letting $E=T_{X} \otimes \mathcal{O}_{X}(K_{X} + nA)$ with $A=2K_X + C_n L$ already yields a bunch of inequalities of Chern classes. This is the starting point of our argument.

\begin{proof}
We first prove the first statement.

By Lemma \ref{veryampleness}, there is some $C_n >0$ such that
\begin{itemize}
  \item $A \triangleq 2K_{X} + C_{n}L$ is very ample.
\end{itemize}
Applying the very ampleness of $A$ in Lemma \ref{tangent twist} shows that
\begin{itemize}
  \item $T_{X} \otimes \mathcal{O}_{X}(K_{X} + nA)$ is nef.
\end{itemize}
Therefore, by Lemma \ref{nef chernclass}, for any partition $\lambda = (\lambda_1,...,\lambda_r)\in P(d,n)$,
\begin{align}\label{start eq}
c_{1}(T_{X} \otimes \mathcal{O}_{X}(K_{X} + nA))^{d} \cdot L^{n-d} \ge c_{\lambda}(T_{X} \otimes \mathcal{O}_{X}(K_{X} + nA)) \cdot L^{n-d} \ge 0 .
\end{align}

By the definition of $A$ and the property of Chern classes,
\begin{align*}
 &c_{\lambda_{i}}(T_{X} \otimes \mathcal{O}_{X}(K_{X} + nA))\\
 &= \sum_{0 \le j \le \lambda_{i}}^{} \binom{n-j}{\lambda_{i}-j} c_{j}(X) \cdot ((2n+1)K_{X} + nC_{n}L)^{\lambda_{i}-j}\\
 &=\sum_{\substack{0 \le j \le \lambda_{i} \\ 0 \le k \le \lambda_{i}-j}} \binom{n-j}{\lambda_{i}-j} \binom{\lambda_{i}-j}{k} (2n+1)^{\lambda_{i}-j-k} n^{k} C_{n}^{k} (-1)^{\lambda_{i}-j-k} c_{j}(X) \cdot c_{1}(X)^{\lambda_{i}-j-k} \cdot L^{k}.
\end{align*}
This shows that $c_{\lambda}(T_{X} \otimes \mathcal{O}_{X}(K_{X} + nA))$ can be written as
\begin{align*}
  c_{\lambda}(T_{X} \otimes \mathcal{O}_{X}(K_{X} + nA)) &= \prod_{i=1}^{r} c_{\lambda_{i}}(T_{X} \otimes \mathcal{O}_{X}(K_{X} + nA))\\
  &= C_{\lambda} c_{\lambda}(X) + \sum_{\mu < \lambda} C_{\lambda,\mu,n}  c_{\mu}(X)\cdot L^{d-|\mu|},
\end{align*}
where $C_{\lambda}, C_{\lambda,\mu,n} \in \mathbb{Z}$ are constants only depending on $n,\lambda,\mu$ and the notion for partitions $\mu < \lambda$ is explained as follows. Note that in the above expression every $c_{\mu}(X)$ is of the form
\begin{equation*}
  c_{\mu}(X) = c_{j_1}(X) \cdot...\cdot c_{j_r}(X) c_{1}(X)^{d-(\sum_{i=1} ^r j_i) - (\sum_{i=1} ^r k_i)},
\end{equation*}
where $0\leq j_i \leq \lambda_i$ for any $i$, $\sum_{i=1} ^r j_i <d$ and $0\leq k_i \leq \lambda_i - j_i $. Therefore, every $\mu$ is obtained by a finite sequence of operators acting on $\lambda$:
\begin{align*}
&(a_{1},...,a_{l}) \longmapsto (+1,a_{1},...,a_{i-1},a_{i}-1,a_{i+1},...,a_{l}),\\
&(a_{1},...,a_{l}) \longmapsto (+0, a_{1},...,a_{j-1},a_{j}-1,a_{j+1},...,a_{l}),
\end{align*}
where $+1, +0$ mean adding one or zero factor $c_1 (X)$ to the Chern monomial.
Given a partition $\nu = (\nu_1,...,\nu_l)$ with every $\nu_i \geq 1$ (or equivalently, a monomial Chern class $c_\nu (X)$ with every $\nu_i \geq 1$), we measure its distance with the power of $c_1 (X)$ by
\begin{equation*}
  d_\nu = \sum_{i=1} ^ l (\nu_i -1).
\end{equation*}
Then it is clear that in the above expansion of $c_{\lambda}(T_{X} \otimes \mathcal{O}_{X}(K_{X} + nA))$, $d_\mu < d_\lambda$. This is what the notion $\mu < \lambda$ means.

Note that
\begin{align*}
 c_{1}(T_{X} \otimes \mathcal{O}_{X}(K_{X} + nA)) = n((2n+1)K_{X} + nC_{n}L) - K_{X},
\end{align*}
Then (\ref{start eq}) can be rewritten as
\begin{align}\label{start eq2}
\left [ n((2n+1)K_{X} + nC_{n}L) - K_{X} \right ]^{d} \cdot L^{n-d} \ge C_{\lambda} c_{\lambda}(X) \cdot L^{n-d} + \sum_{\mu < \lambda} C_{\lambda,\mu,n} c_{\mu}(X) \cdot L^{n- \left | \mu \right |} \ge 0.
\end{align}

Next we apply the induction on $d_\lambda$. When $c_\lambda (X)$ is given by the power of $c_1 (X)$, i.e., $d_\lambda =0$, the result follows by setting $$P_{\lambda}^{\pm}(L^{n},K_{X} \cdot L^{n-1})= \pm (-K_{X})^{d} \cdot L^{n-d}.$$
By induction, we can assume that for every $\mu$ in (\ref{start eq2}), we have linear polynomials $P_{\mu}^{\pm}(x_{0},...,x_{\left | \mu \right |})$ satisfying that
\[
P_{\mu}^{-}(K_{X}^{i} \cdot L^{n-i} : 0 \le i \le \left | \mu \right |) \le c_{\mu}(X) \cdot L^{n-\left | \mu \right |} \le P_{\mu}^{+}(K_{X}^{i} \cdot L^{n-i} : 0 \le i \le \left | \mu \right |).
\]
Regarding the signs of $C_{\lambda,\mu,n}$ in (\ref{start eq2}), we define
$$P_{\lambda,\mu,1}(K_{X}^{i} \cdot L^{n-i} : 0 \le i \le \left | \mu \right |) \triangleq
\begin{cases}
P_{\mu}^{-}(K_{X}^{i} \cdot L^{n-i} : 0 \le i \le \left | \mu \right |), & C_{\lambda,\mu,n} \ge 0,\\
P_{\mu}^{+}(K_{X}^{i} \cdot L^{n-i} : 0 \le i \le \left | \mu \right |), & C_{\lambda,\mu,n} < 0,
\end{cases}$$
and
$$P_{\lambda,\mu,2}(K_{X}^{i} \cdot L^{n-i} : 0 \le i \le \left | \mu \right |) \triangleq
\begin{cases}
P_{\mu}^{+}(K_{X}^{i} \cdot L^{n-i} : 0 \le i \le \left | \mu \right |), & C_{\lambda,\mu,n} \ge 0\\
P_{\mu}^{-}(K_{X}^{i} \cdot L^{n-i} : 0 \le i \le \left | \mu \right |), & C_{\lambda,\mu,n} < 0,
\end{cases}$$
then we have the estimates as follows.
\begin{description}
  \item[Case 1] $C_{\lambda} > 0$.
\end{description}
In this case, by the first inequality in (\ref{start eq2}) the upper bound is given by
\begin{align*}
&c_{\lambda}(X) \cdot L^{n-d}\\
 &\le (C_{\lambda})^{-1} \left [ n((2n+1)K_{X} + nC_{n}L) - K_{X} \right ]^{d} \cdot L^{n-d}
- (C_{\lambda})^{-1} \sum_{\mu < \lambda} C_{\lambda,\mu,n} c_{\mu}(X) \cdot L^{n- \left | \mu \right |} \\
&\le (C_{\lambda})^{-1} \left [ n((2n+1)K_{X} + nC_{n}L) - K_{X} \right ]^{d} \cdot L^{n-d} - (C_{\lambda})^{-1} \sum_{\mu < \lambda} C_{\lambda,\mu,n} P_{\lambda,\mu,1}(K_{X}^{i} \cdot L^{n-i} : 0 \le i \le \left | \mu \right |).
\end{align*}
Note that the above upper bound is a linear polynomial of $K_{X}^{i} \cdot L^{n-i}, 0 \le i \le d$, which we denote by $P_{\lambda}^{+}(K_{X}^{i} \cdot L^{n-i} : 0 \le i \le d)$.

By the second inequality in (\ref{start eq2}), the lower bound is given by
\begin{align*}
&c_{\lambda}(X) \cdot L^{n-d}\\
 &\ge - (C_{\lambda})^{-1} \sum_{\mu < \lambda} C_{\lambda,\mu,n} c_{\mu}(X) \cdot L^{n- \left | \mu \right |}\\
&\ge - (C_{\lambda})^{-1} \sum_{\mu < \lambda} C_{\lambda,\mu,n} P_{\lambda,\mu,2}(K_{X}^{i} \cdot L^{n-i} : 0 \le i \le \left | \mu \right |),
\end{align*}
which gives us the desired bound $P_{\lambda}^{-}(K_{X}^{i} \cdot L^{n-i} : 0 \le i \le d)$.

\begin{description}
  \item[Case 2] $C_{\lambda} < 0$.
\end{description}
This case is similar, we have the upper bound as
\begin{align*}
&c_{\lambda}(X) \cdot L^{n-d}\\
 &\le - (C_{\lambda})^{-1} \sum_{\mu < \lambda} C_{\lambda,\mu,n} c_{\mu}(X) \cdot L^{n- \left | \mu \right |}  \\
&\le - (C_{\lambda})^{-1} \sum_{\mu < \lambda} C_{\lambda,\mu,n} P_{\lambda,\mu,2}(K_{X}^{i} \cdot L^{n-i} : 0 \le i \le \left | \mu \right |) \\
&= P_{\lambda}^{+}(K_{X}^{i} \cdot L^{n-i} : 0 \le i \le d).
\end{align*}
and the lower bound given by
\begin{align*}
&c_{\lambda}(X) \cdot L^{n-d}\\
&\ge (C_{\lambda})^{-1} \left [ n((2n+1)K_{X} + nC_{n}L) - K_{X} \right ]^{d} \cdot L^{n-d} -  (C_{\lambda})^{-1} \sum_{\mu < \lambda} C_{\lambda,\mu,n} c_{\mu}(X) \cdot L^{n- \left | \mu \right |} \\
&\ge (C_{\lambda})^{-1} \left [ n((2n+1)K_{X} + nC_{n}L) - K_{X} \right ]^{d} \cdot L^{n-d} - (C_{\lambda})^{-1} \sum_{\mu < \lambda} C_{\lambda,\mu,n} P_{\lambda,\mu,1}(K_{X}^{i} \cdot L^{n-i} : 0 \le i \le \left | \mu \right |) \\
&= P_{\lambda}^{-}(K_{X}^{i} \cdot L^{n-i} : 0 \le i \le d).
\end{align*}

In any case, we obtain the linear polynomials $P_{\lambda}^{\pm}$. This  finishes the proof of the first statement.

Next we give the proof of the second statement.

To obtain the estimates for $K_{X}^{i} \cdot L^{n-i}$, $1 \le i \le n$, we shall argue by induction on $i$. When $i=1$, it is trivial. Fix some $2 \le i \le n$, and assume that for all $j<i$, we have polynomials $R_{j}^{\pm}(x,y)$ satisfying that $\deg R_{j}^{\pm} (x,y) \le j$ and
\begin{equation}\label{rj}
R_{j}^{-}(L^{n},K_{X} \cdot L^{n-1}) \le (K_{X}^{j} \cdot L^{n-j}) (L^{n})^{j-1} \le R_{j}^{+}(L^{n},K_{X} \cdot L^{n-1}).
\end{equation}

Note that by Lemma \ref{veryampleness}, $K_{X} + (n+1)L$ is nef. Applying Khovanskii-Teissier inequalities implies that
\begin{equation}\label{kt ineq}
0 \le (K_{X} + (n+1)L)^{i} \cdot L^{n-i} \le \frac{( (K_{X} + (n+1)L) \cdot L^{n-1} )^{i}}{(L^{n})^{i-1}}.
\end{equation}
Write
\[
(K_{X} + (n+1)L)^{i} \cdot L^{n-i} = K_{X}^i \cdot L^{n-i} + \sum_{j=0}^{i-1} \binom{i}{j}(n+1)^{i-j}K_{X}^{j} \cdot L^{n-j}.
\]
Then by (\ref{rj}) we have the lower bound given by
\begin{align*}
&K_{X}^{i} \cdot L^{n-i}\\
& \ge - (n+1)^{i}L^{n} - \sum_{j=1}^{i-1} \binom{i}{j}(n+1)^{i-j} K_{X}^{j} \cdot L^{n-j}  \\
&\ge -(n+1)^{i}L^{n} - \sum_{j=1}^{i-1} \binom{i}{j}(n+1)^{i-j} \frac{R_{j}^{+}(L^{n},K_{X} \cdot L^{n-1})}{(L^{n})^{j-1}} \\
&= \frac{R_{i}^{-}(L^{n},K_{X} \cdot L^{n-1})}{(L^{n})^{i-1}},
\end{align*}
where in the last line $R_{i}^{-}(L^{n},K_{X} \cdot L^{n-1})$ is the polynomial given by the polynomials $R_{j}^{+}, j<i$ and $L^n$.
By (\ref{rj}) and (\ref{kt ineq}),  the upper bound is given by
\begin{align*}
&K_{X}^{i} \cdot L^{n-1}\\
& \le \frac{((n+1)L^{n} + K_{X} \cdot L^{n-1})^{i}}{(L^{n})^{i-1}} - (n+1)^{i}L^{n} - \sum_{j=1}^{i-1} \binom{i}{j}(n+1)^{i-j} K_{X}^{j} \cdot L^{n-j}  \\
&\le \frac{((n+1)L^{n} + K_{X} \cdot L^{n-1})^{i}}{(L^{n})^{i-1}} - (n+1)^{i}L^{n} - \sum_{j=1}^{i-1} \binom{i}{j}(n+1)^{i-j} \frac{R_{j}^{-}(L^{n},K_{X} \cdot L^{n-1})}{(L^{n})^{j-1}} \\
&= \frac{R_{i}^{+}(L^{n},K_{X} \cdot L^{n-1})}{(L^{n})^{i-1}},
\end{align*}
where $R_{i}^{+}(L^{n},K_{X} \cdot L^{n-1})$ is the polynomial given by the polynomials $R_{j}^{-}, j<i$ and $L^n, K_{X} \cdot L^{n-1}$.

Thus we obtain the polynomials $R_{i}^{\pm}(x,y)$ as desired, finishing the proof.

\end{proof}

Now we can prove Theorem \ref{main thrm1}.

\begin{proof}[Proof of Theorem \ref{main thrm1}]
Assume that the linear polynomials in Theorem \ref{thrm linear} are given by
$$P_{\lambda}^{\pm}(x_{0},...,x_{d}) = \sum_{i=0}^{d} C_{\lambda,n,i}^{\pm} x_i.$$

We define for $i\ge 1$,
$$R_{\lambda,i,1} \triangleq
\begin{cases}
R_{i}^{+}(L^{n},K_{X} \cdot L^{n-1}), & C_{\lambda,n,i}^{+} \ge 0, \\
R_{i}^{-}(L^{n},K_{X} \cdot L^{n-1}), & C_{\lambda,n,i}^{+} < 0,
\end{cases}
$$
and
$$R_{\lambda,i,2} \triangleq
\begin{cases}
R_{i}^{-}(L^{n},K_{X} \cdot L^{n-1}), & C_{\lambda,n,i}^{-} \ge 0, \\
R_{i}^{+}(L^{n},K_{X} \cdot L^{n-1}), & C_{\lambda,n,i}^{-} < 0.
\end{cases}
$$

Then by Theorem \ref{thrm linear} we have the upper bound
\begin{align*}
&c_{\lambda}(X) \cdot L^{n - d} \le \sum_{i=0}^{d} C_{\lambda,n,i}^{+} K_{X}^{i} \cdot L^{n-i}  \\
&\le C_{\lambda,n,0}^{+} L^n + \sum_{i=1}^{d} C_{\lambda,n,i}^{+} \frac{R_{\lambda,i,1}(L^{n},K_{X} \cdot L^{n-1})}{(L^n)^{i-1}} \\
&= \frac{Q_{\lambda}^{+}(L^{n},K_{X} \cdot L^{n-1})}{(L^{n})^{d- 1}},
\end{align*}
where $Q_{\lambda}^{+}(L^{n},K_{X} \cdot L^{n-1})$ is the polynomial given by the polynomials $R_{\lambda,i,1}$ and $L^n$,
and we have the lower bound
\begin{align*}
&c_{\lambda}(X) \cdot L^{n - d} \ge \sum_{i=0}^{d} C_{\lambda,n,i}^{-} K_{X}^{i} \cdot L^{n-i}  \\
&\ge C_{\lambda,n,0}^{-} L^n +\sum_{i=1}^{d} C_{\lambda,n,i}^{-} \frac{R_{\lambda,i,2}(L^{n},K_{X} \cdot L^{n-1})}{(L^n)^{i-1}} \\
&= \frac{Q_{\lambda}^{-}(L^{n},K_{X} \cdot L^{n-1})}{(L^{n})^{d - 1}},
\end{align*}
where $Q_{\lambda}^{-}(L^{n},K_{X} \cdot L^{n-1})$ is the polynomial given by the polynomials $R_{\lambda,i,2}$ and $L^n$.

Therefore we obtain polynomials $Q_{\lambda}^{\pm}(x,y)$ satisfying that $\deg Q_{\lambda}^{\pm}(x,y) \le d$ and
\begin{align*}
Q_{\lambda}^{-}(L^{n},K_{X} \cdot L^{n-1}) \le (c_{\lambda}(X) \cdot L^{n - d}) (L^{n})^{d - 1} \le Q_{\lambda}^{+}(L^{n},K_{X} \cdot L^{n-1}).
\end{align*}
Assume that $$Q_{\lambda}^{\pm}(x,y) =  \sum_{\substack{i+j \le d \\ i,j \ge 0 }} a_{\lambda,ij}^{\pm} x^{i} y^{j},$$ and define
\begin{align*}
Q_{\lambda}(x,y) \triangleq \sum_{\substack{i+j \le d \\ i,j \ge 0 }} a_{\lambda,ij} x^{i} y^{j} + \sum_{\substack{i+j \le d \\ i,j \ge 0 }} \left | a_{\lambda,ij}^{+} + a_{\lambda,ij}^{-} \right | (n+1)^{j} x^{i+j},
\end{align*}
where $$a_{\lambda,ij} \triangleq \max \left \{ a_{\lambda,ij}^{+}, - a_{\lambda,ij}^{-} \right \}.$$

We claim that $$\max \left \{ Q_{\lambda}^{+}(L^{n},K_{X} \cdot L^{n-1}), - Q_{\lambda}^{-}(L^{n},K_{X} \cdot L^{n-1}) \right \} \le Q_{\lambda}(L^{n},K_{X} \cdot L^{n-1}).$$

To this end, note that this is obviously true when $K_{X} \cdot L^{n-1} \ge 0$.
We only need to consider the case when $K_{X} \cdot L^{n-1} < 0$. For $j$ even,
\begin{align*}
\max \left \{ a_{\lambda,ij}^{+} (L^{n})^{i} (K_{X} \cdot L^{n-1})^{j}, - a_{\lambda,ij}^{-} (L^{n})^{i} (K_{X} \cdot L^{n-1})^{j} \right \} \le a_{\lambda,ij} (L^{n})^{i} (K_{X} \cdot L^{n-1})^{j}.
\end{align*}
For $j$ odd,
\begin{align*}
&\max \left \{ a_{\lambda,ij}^{+} (L^{n})^{i} (K_{X} \cdot L^{n-1})^{j}, - a_{\lambda,ij}^{-} (L^{n})^{i} (K_{X} \cdot L^{n-1})^{j} \right \}\\
& = \min \left \{ a_{\lambda,ij}^{+} , - a_{\lambda,ij}^{-} \right \} (L^{n})^{i} (K_{X} \cdot L^{n-1})^{j}  \\
&= \left ( \min \left \{ a_{\lambda,ij}^{+}, - a_{\lambda,ij}^{-} \right \} + \left | a_{\lambda,ij}^{+} + a_{\lambda,ij}^{-} \right | \right ) (L^{n})^{i} (K_{X} \cdot L^{n-1})^{j} - \left | a_{\lambda,ij}^{+} + a_{\lambda,ij}^{-} \right | (L^{n})^{i} (K_{X} \cdot L^{n-1})^{j} \\
&= a_{\lambda,ij} (L^{n})^{i} (K_{X} \cdot L^{n-1})^{j} - \left | a_{\lambda,ij}^{+} + a_{\lambda,ij}^{-} \right | (L^{n})^{i} (K_{X} \cdot L^{n-1})^{j} \\
&\le a_{\lambda,ij} (L^{n})^{i} (K_{X} \cdot L^{n-1})^{j} - \left | a_{\lambda,ij}^{+} + a_{\lambda,ij}^{-} \right | (L^{n})^{i} (-(n+1)L^{n})^{j} \\
&= a_{\lambda,ij} (L^{n})^{i} (K_{X} \cdot L^{n-1})^{j} + \left | a_{\lambda,ij}^{+} + a_{\lambda,ij}^{-} \right | (n+1)^{j} (L^{n})^{i+j},
\end{align*}
where in the second equality we apply the identity
$$\max(a, -b)=|a+b|+\min(a, -b),$$
and in the fourth inequality we
use the fact that $K_{X} + (n+1)L$ is nef, implying that
\begin{align*}
(K_{X} + (n+1)L) \cdot L^{n-1} =K_{X} \cdot L^{n-1} + (n+1)L^{n} \ge 0.
\end{align*}
This finishes the proof of our claim.

Therefore, we have
\[
\left | c_{\lambda}(X) \cdot L^{n- d } \right | \le \frac{\max \left \{ Q_{\lambda}^{+}(L^{n},K_{X} \cdot L^{n-1}), - Q_{\lambda}^{-}(L^{n},K_{X} \cdot L^{n-1}) \right \}}{(L^{n})^{d - 1}} \le \frac{Q_{\lambda}(L^{n},K_{X} \cdot L^{n-1})}{(L^{n})^{d - 1}}.
\]

\end{proof}

\section{Applications and extensions}

\subsection{Boundedness for Chern numbers}
In the recent work by Du-Sun \cite{durongChernIneq}, by using the method of pulling back Schubert classes in the Chow group of a Grassmannian under the Gauss map, it was proved that: for a projective manifold $X$ over an algebraically closed field $k$, if $K_X$ or $-K_X$ is ample and $k$ is of characteristic 0, or if $K_X$ or $-K_X$ is ample and globally generated and $k$ is of positive characteristic, then there exists a constant $c_n$ depending only on $n$ such that for any monomial Chern classes of top degree, the Chern number ratios satisfies the following inequality
\begin{equation*}
\left|\frac{c_\lambda(X)}{c_1 (X) ^{n}}\right|\leq c_n.
\end{equation*}

Theorem \ref{main thrm1} generalizes their result over fields of characteristic 0. In particular, letting $L=K_X$ or $L=-K_X$ immediately recovers \cite{durongChernIneq}.

\begin{cor}
Let $X$ be a projective manifold of dimension $n$, with $K_{X}$ or $-K_{X}$ ample, then the Chern number ratios
$$\left [ \frac{c_{\lambda}(X)}{c_{1}(X)^{n}} \right ]_{\lambda \in P(n,n)} \in \mathbb{R}^{p(n)}$$
is contained in a bounded set in $\mathbb{R}^{p(n)}$ independent of $X$, where $p(n)$ is the partition number.
\end{cor}

\begin{proof}
We only consider the case when $K_X$ is ample. The case when $-K_{X}$ is ample is similar.

For $\lambda \in P(n,n)$, write $Q_\lambda (x, x)=\sum_{i=0}^n b_i x^i$. Letting $L=K_X$ in Theorem \ref{main thrm1} and using $L^n \geq 1$ imply that
\begin{align*}
\left| \frac{c_\lambda (X)}{c_1 (X) ^n}\right|\leq \frac{Q_\lambda (L^n, L^n)}{(L^n)^n}
=\sum_{i=0}^n b_i \left(\frac{1}{L^n}\right)^{n-i}
\leq \sum_{i=0}^n |b_i|.
\end{align*}
Note that the $b_i$ depend only on $n, \lambda$ and there is only finite number (depending only on $n$) of partitions $\lambda$. This finishes the proof.
\end{proof}

Indeed, we have:
\begin{cor}\label{upperbd chern}
There is some uniform constant $c(n, v, w) >0$ such that for any projective manifold $X$ of dimension $n$ and any ample line bundle $L$ on $X$ with
\begin{itemize}
  \item $L^n \leq v$,
  \item $K_X \cdot L^{n-1} \leq w$,
\end{itemize}
we have
\begin{equation*}
  \left|c_\lambda(X)\cdot L^{n -d}\right|\leq c(n, v, w).
\end{equation*}

\end{cor}

\begin{proof}
By Lemma \ref{veryampleness}, $L^n \leq v$ implies that
\begin{equation*}
  K_X \cdot L^{n-1} \geq -(n+1)v,
\end{equation*}
thus $-(n+1)v \leq K_X \cdot L^{n-1} \leq w$.
Now the result follows from Theorem \ref{main thrm1}.
\end{proof}

\subsection{Asymptotic Riemann-Roch type inequalities}
Another consequence is an asymptotic version of the sharper Riemann-Roch type inequality for ample line bundles.

\begin{cor}
There exists a polynomial $Q(z)$ with $\deg Q \le n-m-1$, such that for any projective manifold $X$ of dimension $n$ and any ample line bundle $L$ on $X$, for $k$ large enough,
$$\left| h^{0}(X,kL)-\sum_{i=0}^m a_i k^{n-i}\right | \le Q(k),$$
where $\sum_{i=0}^m a_i k^{n-i}$ is the truncation of the Hilbert polynomial
$$ \chi(X,kL) = \int_X ch(kL) \cdot td(X)=\sum_{i=0}^n a_i k^{n-i}.$$
The coefficients of $Q(z)$ are determined by $L^{n},K_{X} \cdot L^{n-1}, n$ and the partitions $\lambda$ with $|\lambda|\geq m+1$.
\end{cor}

\begin{proof}
For sufficiently large $r$, $h^{0}(X,rL) = \chi(rL)$.
Write $$ch(kL) = \sum_{i=0}^{n} \frac{k^{i}}{i!} L^{i},\ td(X) = \sum_{i=0}^{n} S_{i}(X).$$
Note that every $S_{i}(X)$ is given by some combination of $c_{\lambda}(X)$ with $| \lambda | = i$. Therefore,
\begin{align*}
a_i = \frac{1}{(n - i)!} S_{i}(X) \cdot L^{n-i}.
\end{align*}

Then the result follows from Theorem \ref{main thrm1} immediately.
\end{proof}

\subsection{Logarithmic tangent bundles}

In this section, using similar method, we establish the logarithmic analogy of Theorem \ref{main thrm1}.
Let $X$ be a projective manifold of dimension $n$, and let $D=\sum_{i=1} ^N D_i$ be a reduced SNC divisor on $X$. Denote by $\Omega_X ^1 (\log D)$ the logarithmic cotangent bundle of the pair $(X, D)$. Note that we have the following exact sequence:
\begin{align*}
  &0\rightarrow \Omega_X ^1 \rightarrow \Omega_X ^1 (\log D) \rightarrow \oplus_{i=1} ^N \mathcal{O}_{D_i} \rightarrow 0, \\
\end{align*}
where $\mathcal{O}_{D_i}$ is the sheaf given by
\begin{align*}
  &0\rightarrow \mathcal{O}_{X} (-D_i) \rightarrow \mathcal{O}_{X} \rightarrow  \mathcal{O}_{D_i} \rightarrow 0. \\
\end{align*}
Then we have that
\begin{equation}\label{logchern}
\begin{aligned}
  c_k (\Omega_X ^1 (\log D)) &= \sum_{i=0} ^k c_{k-i} (\Omega_X ^1)\cdot \left( \sum_{j_1 +...+j_N =i} D_1 ^{j_1}\cdot...\cdot  D_N ^{j_N}\right)\\
  &=\sum_{i=0} ^k c_{k-i} (\Omega_X ^1) \cdot \sum_{|J|=i} D^J,
\end{aligned}
\end{equation}
where for the index $J=(j_1,...,j_N)\in \mathbb{N}^N$ we denote $D^J = D_1 ^{j_1}\cdot...\cdot  D_N ^{j_N}$ and $|J|= j_1 +...+j_N$.

We intend to prove:

\begin{thrm}\label{thrm log}
Let $X$ be a projective manifold of dimension $n$, and let $D=\sum_{i=1} ^N D_i$ be a reduced SNC divisor on $X$. Denote by $\Omega_X ^1 (\log D)$ the logarithmic cotangent bundle of the pair $(X, D)$. Let $L$ be an ample line bundle on $X$.
Then for any partition $\lambda \in P(d, n)$, the Chern number $c_\lambda (\Omega_X ^1 (\log D))\cdot L^{n-d}$ is polynomially controlled by the intersection numbers $$x_J=D^J\cdot L^{n-|J|},\  y_J=K_X \cdot D^J\cdot L^{n-|J|-1},$$
where $J$ varies among the index set $\{J= (j_1,...,j_N): 0\leq j_k \leq 2\}$. More precisely, there is a positive integer $k(\lambda, n)$ and
a family of polynomials $\{Q_{\lambda} ^l\}_{l=1} ^{k(\lambda, n)}$ such that
\begin{align*}
  \min_{1\leq l\leq k(\lambda, n)} \{Q_{\lambda} ^l(x_J, y_J)\}
  \leq c_\lambda (\Omega_X ^1 (\log D))\cdot L^{n-d}
  \leq \max_{1\leq l\leq k(\lambda, n)} \{Q_{\lambda} ^l (x_J, y_J)\}.
\end{align*}

\end{thrm}

\begin{rmk}
When $N=1$, i.e., $D$ is irreducible, the theorem shows that the Chern number $c_\lambda (\Omega_X ^1 (\log D))\cdot L^{n-d}$ is polynomially controlled by the following six quantities $$x_k= D^{k}\cdot L^{n-k},\ y_k= K_{X}\cdot D^{k}\cdot L^{n-1-k}, \ 0\le k\le 2.$$
\end{rmk}

Next we give the proof of Theorem \ref{thrm log}, where the method is rather similar to Theorem \ref{main thrm1}.

Given a partition $\lambda=(\lambda_{1},...,\lambda_{k})\in P(d,n)$ with $1 \le \lambda_{i} \le n$, by (\ref{logchern}) we have
\begin{equation}\label{expanding eq}
\begin{aligned}
    &c_{\lambda} (\Omega_{X}^{1} (\log D)) \cdot L^{n-|\lambda|} = \left [ \prod_{l=1}^{k} ( c_{\lambda_{l}} (\Omega_{X}^{1}) + \sum_{\substack{
    0\le i<\lambda_{l} \\ |I|=\lambda_{l}-i}} c_{i} (\Omega^{1}_{X}) \cdot D^{I} ) \right ] \cdot L^{n-|\lambda|}\\
 & = c_{\lambda} (\Omega^{1}_{X}) \cdot L^{n-\left | \lambda \right |} + \sum_{\substack{
    \mu \prec \lambda \\ |J|=\left | \lambda \right | - \left | \mu \right |} } A^{\lambda}_{\mu,j} c_{\mu} (\Omega^{1}_{X}) \cdot D^{J} \cdot L^{n-\left | \mu \right | -|J|}
\end{aligned}
\end{equation}
where $\mu \prec \lambda$ means that $\mu=(\mu_{1},...,\mu_{k})$ satisfies $0\le\mu_{i}\le\lambda_{i}$ and $\left | \mu \right | < \left | \lambda \right |$.

To give a clean and clear exposition, we first consider the case when $N=1$.

\begin{description}
  \item[Special case] $N=1$.
\end{description}

By Theorem \ref{main thrm1}, we have the bounds for $c_{\lambda} (\Omega^{1}_{X}) \cdot L^{n-\left | \lambda \right |}$. It remains to estimate
$$c_{\mu} (\Omega^{1}_{X}) \cdot D^{j} \cdot L^{n-\left | \mu \right | -j}$$
where $\mu \prec \lambda, j=\left | \lambda \right | - \left | \mu \right |\geq 1$. This can be proved by an analog of Theorem \ref{thrm linear}.

\begin{lem}\label{linearpoly}
Given $j\ge 1, \mu$ as above, there exist linear polynomials $\left\{ P^{\pm}_{\mu,j} (x_{i,l}) : 0\le i\le \left | \mu \right |, 1\le l\le j \right\}$ such that
    \begin{align*}
       P^{-}_{\mu,j} (K_{X}^{i} \cdot D^{l} \cdot L^{n-i-l}: 0\le i\le \left | \mu \right |, 1\le l\le j)   &\le c_{\mu} (\Omega^{1}_{X}) \cdot D^{j} \cdot L^{n-\left | \mu \right | -j} \\
        &\le  P^{+}_{\mu,j} (K_{X}^{i} \cdot D^{l} \cdot L^{n-i-l} : 0\le i\le \left | \mu \right |, 1\le l\le j )
    \end{align*}
\end{lem}

\begin{proof}
Note that $D$ is a smooth hypersurface on $X$, therefore $(K_{X}+D)|_{D} = K_{D}$. By the results in Section \ref{prel},
\begin{itemize}
  \item $K_{D}+(n+1)L_{|D} = (K_{X}+D+(n+1)L)_{|D}$ is ample on $D$,
  \item $(T_{X} \otimes \mathcal{O} ((2n+1)K_{X}+nC_{n}L))_{|D}$ is nef on $D$, where $C_{n}$ is the constant as in Lemma \ref{veryampleness}.
\end{itemize}

Applying Lemma \ref{nef chernclass} implies that
\begin{equation}\label{fllog}
\begin{aligned}
    0 &\le c_{\mu} (T_{{X}_{|D}} \otimes \mathcal{O}_{X}((2n+1)K_{X}+nC_{n}L)_{|D}) \cdot (D+K_{X}+(n+1)L)_{|D}^{j-1} \cdot L_{|D}^{n-\left | \mu \right |-j}\\
&\le c_{1} (T_{{X}_{|D}} \otimes \mathcal{O}_{X}((2n+1)K_{X}+nC_{n}L)_{|D})^{\left | \mu \right |} \cdot (D+K_{X}+(n+1)L)_{|D}^{j-1} \cdot L_{|D}^{n-\left | \mu \right |-j}.
\end{aligned}
\end{equation}

Similarly to rewriting (\ref{start eq}) in the form of (\ref{start eq2}),
(\ref{fllog}) can be expanded as follows:
\begin{equation}\label{relative fl}
\begin{aligned}
    0 & \le A_{\mu,j} c_{\mu}(\Omega^{1}_{X}) \cdot D^{j} \cdot L^{n-\left | \mu \right | -j} + \sum_{\rho< \mu} C^{\mu,j}_{\rho} c_{\rho}(\Omega^{1}_{X}) \cdot D^{j} \cdot L^{n-\left | \rho \right | -j}\\
    &+ \sum_{\substack{
    1\le l< j \\ \left | \rho \right | \le \left | \mu \right | +j-l
    }} C^{\mu,l}_{\rho}c_{\rho}(\Omega^{1}_{X}) \cdot D^{l} \cdot L^{n-\left | \rho \right |-l}\\
&\triangleq(I) + (II) +(III)\\
    &\leq \sum_{0\le i\le \left | \mu \right |} A^{\mu,j}_{i} K_{X}^{i} \cdot D^{j}\cdot L^{n-i-j} + \sum_{\substack{
    1\le l< j \\ 0\le i\le \left | \mu \right | +j-l
    }} A^{\mu,l}_{i} K_{X}^{i}\cdot D^{l}\cdot L^{n-i-l}.
\end{aligned}
\end{equation}

Then a double induction on $j$ and $\rho< \mu$ implies the desired estimates for the terms $(II), (III)$, which in turn gives the bound for the intersection number $c_{\mu}(\Omega^{1}_{X}) \cdot D^{j} \cdot L^{n-\left | \mu \right | -j}$.

This finishes the proof of the lemma.

\end{proof}

Given $i, j$ such that $1\leq j \leq n-i$, we show that $K_{X}^{i}\cdot D^{j}\cdot L^{n-i-j}$ can be polynomially controlled by
\begin{equation*}
  x_k= D^{k}\cdot L^{n-k},\ y_k= K_{X}\cdot D^{k}\cdot L^{n-1-k}
\end{equation*}
where $1\leq k \leq 2$.

\begin{lem}\label{rpoly}
    For any given $i, j$ satisfying that $1\leq j \leq n-i$, there exist a positive integer $k(i,j,n)$ and a family of polynomials $\{ R_{i,j}^{l} (x_{k},y_{k}):k=1,2 \}_{l=1} ^{ k(i,j,n)} $ such that,
\begin{align*}
        \min_{1\leq l\le k(i,j,n)} \{ R_{i,j}^{l} (x_k, y_k:k=1,2) \} \le K_{X}^{i}\cdot D^{j}\cdot L^{n-i-j}
        \le \max_{1\leq l\le k(i,j,n)} \{ R_{i,j}^{l} (x_k, y_k:k=1,2) \}
    \end{align*}
\end{lem}

For $j=1$, from the proof below we will see that the bounds do not depend on the variables $x_2, y_2$.

\begin{proof}
We first consider the case when $j=1$. Using the fact that
\begin{equation}\label{amplebd1}
  L_{|D} ^{n-1} \geq 1
\end{equation}
and the Khovanskii-Teissier inequalities (Lemma \ref{ktineq}), for $i'\ge 1$, we obtain:
\begin{align*}
    0&\le (K_{X}+(n+1)L)_{|D}^{i'} \cdot L_{|D}^{n-i'-1} \\
&\le \left[(K_{X}+(n+1)L)_{|D}^{i'}\cdot L_{|D}^{n-i'-1}\right] (L_{|D}^{n-1})^{i'-1} \\
    &\le \left[(K_{X}+(n+1)L)_{|D}\cdot L_{|D}^{n-2}\right]^{i'},
\end{align*}
which is equivalent to that
\begin{align*}
    0 &\le K_{X}^{i'}\cdot D\cdot L^{n-i'-1} + \sum_{0\le l\le i'-1} a_l K_{X}^{l}\cdot D\cdot L^{n-l-1} \\
&\triangleq (A) +(B)\\
    &\le \sum_{0\le l\le i'} b_l (K_{X}\cdot D\cdot L^{n-2})^{l} (D\cdot L^{n-1})^{i'-l},
\end{align*}
where $a_l =\binom{i'}{l} (n+1)^{i'-l}$ and $b_l=\binom{i'}{l} (n+1)^{i'-l}$.
Applying an induction on $l$ shows that the term $(B)$ satisfies the desired estimate, which then implies the bound for  $K_{X}^{i'}\cdot D\cdot L^{n-i'-1}$. This finishes the proof for $j=1$. Indeed, we obtain better results in this case, that is, there exist polynomials $R_{i',1}^{\pm} (x_1, y_1)$, such that
\begin{align*}
    R_{i',1}^{-} (x_1, y_1) &\le K_{X}^{i'}\cdot D\cdot L^{n-i'-1} \le  R_{i',1}^{+} (x_1, y_1).
\end{align*}

The result for $j\geq 2$ follows from the following claim and a similar inductive process.
\begin{claim}
\begin{equation}\label{ktlog1}
\begin{aligned}
    0 &\le (K_{X}+(n+2)L)_{|D}^{i}\cdot L_{|D}^{n-i-j}\cdot (D+K_{X}+(n+1)L)_{|D}^{j-1}\\
&\le \left[(D+K_{X}+(n+1)L)_{|D}\cdot L_{|D}^{n-3}\cdot (K_{X}+(n+2)L) _{|D}\right] ^{i(j-1)}.
\end{aligned}
\end{equation}
\end{claim}

To prove the claim, first note that $(D+K_{X}+(n+1)L)_{|D}$ is ample on $D$, thus
\begin{equation}\label{amplegeq1}
  L_{|D}^{n-j}\cdot (D+K_{X}+(n+1)L)_{|D}^{j-1} \geq 1,
\end{equation}
which implies that $ (K_{X}+(n+2)L)_{|D}^{i}\cdot L_{|D}^{n-i-j}\cdot (D+K_{X}+(n+1)L)_{|D}^{j-1}$ is bounded above by
\begin{equation}\label{trivbd}
 \left[(K_{X}+(n+2)L)_{|D}^{i}\cdot L_{|D}^{n-i-j}\cdot (D+K_{X}+(n+1)L)_{|D}^{j-1}\right] \left[L_{|D}^{n-j}\cdot (D+K_{X}+(n+1)L)_{|D}^{j-1}\right]^{i-1}.
\end{equation}
An application of Lemma \ref{ktineq} yields that the product (\ref{trivbd}) is bounded above by
\begin{equation}\label{trivbd1}
\left[(K_{X}+(n+2)L)_{|D}\cdot L_{|D}^{n-1-j}\cdot (D+K_{X}+(n+1)L)_{|D}^{j-1}\right]^{i}.
\end{equation}
Multiplying (\ref{trivbd1}) by
\begin{equation}\label{amplebd10}
 [L_{|D} ^{n-2} \cdot (K_X +(n+2)L)_{|D}]^{i(j-2)}\ge 1
\end{equation}
and applying Lemma \ref{ktineq} to the product again, we obtain the claimed estimate in (\ref{ktlog1}).

The inequalities in (\ref{ktlog1}) can be rewritten as
\begin{align*}
    0 &\le K_{X}^{i}\cdot D^{j}\cdot L^{n-i-j} + \sum_{0\le l<i} A_l K_{X}^{l}\cdot D^{j}\cdot L^{n-l-j} + \sum_{1\le k<j} \sum_{ 0\le l\le i+j-k} C_{l,k} K_{X}^{l}\cdot D^{k}\cdot L^{n-l-k}\\
&\triangleq (I)+(II)+(III)\\
    &\le \left [\sum_{k=1,2} \sum_{l\le 3-k} B_{l,k}  K_{X}^{l}\cdot D^{k}\cdot L^{n-l-k} \right ]^{i(j-1)}.
\end{align*}

Applying inductions on $i, j$ implies that the terms $(II), (III)$ and the term of the type $K_X ^2 \cdot D\cdot L^{n-3}$ in the last line satisfy the desired estimates. Using the fact that for $p\geq 0$, $|a^p|\le \max (r_+ ^p, (-r_- )^p)$ whenever $r_- \le a\le r_+$ and some elementary analysis gives the estimate for $K_{X}^{i}\cdot D^{j}\cdot L^{n-i-j}$, finishing the proof for $j\ge 2$.

\end{proof}

The case $N=1$ then follows from Lemma \ref{linearpoly} and Lemma \ref{rpoly}.

\begin{description}
  \item[General case] $N\geq 2$.
\end{description}

The case for general $N$ can be proved by a rather similar way, we omit the details and just give a sketch below.

By Theorem \ref{main thrm1} and (\ref{expanding eq}), we need to estimate $$c_{\mu} (\Omega_{X}^{1})\cdot D^{J}\cdot L^{n-\left | \mu \right |-|J|}$$
where $\mu \prec \lambda$ and $|J| = |\lambda|-|\mu|\geq 1$.
Although $N$ can be large, $c_{\mu} (\Omega_{X}^{1})\cdot D^{J}\cdot L^{n-\left | \mu \right |-|J|}$ involves at most $n$ components of $D$.
Therefore, in the sequel we can assume $N\le n$ and all $j_l >0$ in $J$. Moreover, we can assume that $Z\triangleq \bigcap_{i=1}^{N} D_i$ is nonempty.

To simplify the notations below, when $J=(1,...,1,j_{M+1},...,j_{N})$ with $j_{l}\ge 2,l=M+1,...,N$, we denote
$$\Omega=(D_{M+1} +K_{X}+(n+1)L)_{|Z}^{j_{M+1} -1}\cdot ... \cdot (D_N +K_{X}+(n+1)L)_{|Z}^{j_N -1},$$
and when $J=(1,...,1)$, denote $\Omega=1$.

Analogy to Lemma \ref{linearpoly}, one can prove the following result:

\begin{lem}\label{lemlinearN}
For $0\le i\le \left | \mu \right |, I=(i_1,...,i_N)$ with $j_l \ge i_l \ge 1$, denote
\begin{equation*}
  x_{i,I}=K_{X}^{i} \cdot D^{I} \cdot L^{n-i-\left | I \right |}.
\end{equation*}
Then there exist linear polynomials $P^{\pm}_{\mu,J} (x_{i,I})$ such that
    \begin{align*}
        P^{-}_{\mu,J} (x_{i, I}) \le c_{\mu} (\Omega^{1}_{X}) \cdot D^{J} \cdot L^{n-\left | \mu \right | -\left | J \right |} \le  P^{+}_{\mu,J} (x_{i, I}).
    \end{align*}
\end{lem}

This follows from an inductive process and the analogy to (\ref{fllog}):
\begin{align*}
    0 &\le c_{\mu} ({T_{X}}_{|Z} \otimes \mathcal{O}_{X}((2n+1)K_{X}+nC_{n}L)_{|Z}) \cdot \Omega\cdot L_{|Z}^{n-\left | \mu \right |-\left | J \right |} \\
    &\le c_{1} ({T_{X}}_{|Z} \otimes \mathcal{O}_{X}((2n+1)K_{X}+nC_{n}L)_{|Z})^{\left | \mu \right |} \cdot \Omega \cdot L_{|Z}^{n-\left | \mu \right |-\left | J \right |},
\end{align*}

Similar to Lemma \ref{rpoly}, we have:

\begin{lem}\label{lemmarpolyN}
Fix $i\ge 0, I=(i_1,...,i_N)$ with $i_l \ge 1$ and $i+|I|\leq n$. Denote
\begin{equation*}
  x_J=D^{J}\cdot L^{n-\left | J \right |},\ y_J=K_{X}\cdot D^{J}\cdot L^{n-1-\left | J \right |},
\end{equation*}
where $J$ varies in the set $\{J=(j_1,...,j_N): 1\le j_k \le 2\}$.
Then there exist a positive integer $k(i,I,n)$ and a family of polynomials $\{ R_{i,I}^{l} (x_{J},y_{J})\}_{l=1} ^{k(i, I, n)}$ such that
    \begin{align*}
        \min_{1\le l\le k(i,I,n)} \{ R_{i,I}^{l} (x_J, y_J) \} \le K_{X}^{i}\cdot D^{I}\cdot L^{n-i-\left | I \right |} \le \max_{1\le l\le k(i,I,n)} \{ R_{i,I}^{l} (x_J, y_J) \}.
    \end{align*}

\end{lem}

One can first prove the case when $I=(1,...,1)$ by using
\begin{align}\label{amplebd2}
    0&\le (K_{X}+(n+1)L)_{|Z}^{i} \cdot L_{|Z}^{n-i-N}
    \le \left[(K_{X}+(n+1)L)_{|Z}\cdot L_{|Z}^{n-1-N}\right]^{i}.
\end{align}
For the other case when $I=(1,...,1,i_{M+1},...,i_{N})$ with $i_{l}\ge 2,l=M+1,...,N$, it follows from an inductive process and the analogy to (\ref{ktlog1}):
\begin{equation}\label{amplebd3}
\begin{aligned}
    0 &\le (K_{X}+(n+2)L)_{|Z}^{i}\cdot L_{|Z}^{n-\left | I \right |-i}\cdot \Omega\\
&\le \left [ \left (\prod_{l=M+1}^{N} (D_{l}+K_{X}+(n+1)L)_{|Z} \right ) \cdot L_{|Z}^{n-1-2N+M}\cdot (K_{X}+(n+2)L)_{|Z} \right ] ^{i\prod_{l=M+1}^{N} (i_l -1)},
\end{aligned}
\end{equation}
which can be obtained by an iterative application of Lemma \ref{ktineq}.

Applying Lemmas \ref{lemlinearN} and \ref{lemmarpolyN} implies the result for general $N$.

\begin{rmk}
Regarding the bound for Chern classes of the logarithmic tangent bundle, comparing to the form given by Theorem \ref{main thrm1}, the main difference is about the denominator. The reason is that, in the applications of Khovanskii-Teissier inequalities, we just use the trivial bounds given by the terms like (\ref{amplebd1}), (\ref{amplegeq1}), (\ref{amplebd10}) which considerably simplify the inductive process, and then we obtain the estimates (\ref{ktlog1}), (\ref{amplebd2}) and (\ref{amplebd3}). If we keep track of these terms, it appears unclear to us which kind of exact bounds -- which would have much more advantages -- can be expected.
\end{rmk}

\bibliographystyle{alpha}
\bibliography{RRCreference}

\bigskip

\bigskip

\noindent
\textsc{Tsinghua University, Beijing 100084, China}\\
\noindent
\verb"Email: lux22@mails.tsinghua.edu.cn"\\
\noindent
\verb"Email: jianxiao@tsinghua.edu.cn"

\end{document}